\documentclass{elsart}
\usepackage{graphicx}
\usepackage{graphics}
\usepackage[dvips]{epsfig}

\newcommand{\R}{I \! \! R}
\newcommand{\N}{I \! \! {N}}
\newcommand{\C}{I \! \! \! \! {C}}

\newcommand{\bc}{{\bf c}}

\newcommand{\bud}{{\bf \underline{d}}}

\newcommand{\ba}{{\bf a}}

\newcommand{\brho}{{\mbox{\boldmath $\rho$}}}
\newcommand{\bthe}{{\mbox{\boldmath $\theta$}}}
\newcommand{\bbb}{{\bf b}}
\newcommand{\bux}{{\bf \underline{x}}}
\newcommand{\buc}{{\bf \underline{c}}}
\newcommand{\buy}{{\bf \underline{y}}}

\newcommand{\bV}{{\bf  V}}
\newcommand{\bvs}{{\bf {\underline a}}}

\newcommand{\bxi}{{\mbox{\boldmath $\xi$}}}
\newcommand{\buxi}{{\mbox{\boldmath $\underline{\xi}$}}}
\newcommand{\buzet}{{\mbox{\boldmath $\underline{\zeta}$}}}

\newcommand{\bzet}{{\mbox{\boldmath $\zeta$}}}

\newcommand{\bua}{{\bf{\underline a}}}
\newcommand{\bub}{{\bf{\underline b}}}
\newcommand{\buat}{\tilde{\bf{\underline a}}}

\newcommand{\bue}{{\bf{\underline e}}}
\newcommand{\vs}{{\underline a}}

\newcommand{\umu}{{\underline{\mu}}}

\newcommand{\us}{{\underline s}}
\newcommand{\bU}{{\bf U}}
\newcommand{\bu}{{\bf u}}

\newtheorem{theorem}{Theorem}
\newtheorem{lemma}{Lemma}
\newtheorem{corollary}{Corollary}
\newcommand{\bb}{\begin{eqnarray*}}
\newcommand{\be}{\end{eqnarray*}}

\begin{document}
\begin{frontmatter}

\title{Bivariate one-sample optimal location test for  spherical stable densities by Pade' methods}
\author{P. Barone }
\address{ Istituto per le Applicazioni del Calcolo ''M. Picone'',
C.N.R.,\\
Via dei Taurini 19, 00185 Rome, Italy \\
e-mail: piero.barone@gmail.com, p.barone@iac.cnr.it}

\begin{abstract}

Complex signal detection in additive noise can be performed by a one-sample bivariate location test. Spherical symmetry is assumed for the noise density as well as closedness with respect to linear transformation. Therefore the noise is assumed to have spherical distribution with $\alpha-$stable radial density. In order to cope with this difficult setting the original sample is transformed by Pade' methods giving rise to a new sample with universality properties. The stability assumption is then reduced to the Gaussian one and it is proved that a known van der Waerden type test, with optimal properties, based on the new sample can be used. Furthermore a new test in the same class of optimal tests is proposed which is more powerful that the van der Waerden type one.

\end{abstract}

\begin{keyword}

 random matrices;  Pad\'{e}  approximants; multivariate location
\end{keyword}

\end{frontmatter}


\section*{Introduction}

Additive noise filtering is a common problem in many experimental situations. Sometimes happens however that what matters is to understand if a signal is present or not in the observations. The specific shape of the signal is not relevant. Moreover sometimes it is not possible to make assumptions on the statistical distribution of the noise. The problem consists then in characterizing the noise w.r. to any possible signal with the only constraint that the noise is additive. In the following we assume that the noise can be represented by a discrete time, complex valued,
 stationary process such that every finite set of random variables of the process have a multivariate spherical distribution centered in zero. In the limit case in which this process reduces to a single complex random variable this is equivalent to consider a couple of real random variables with bivariate spherical distribution.  As the noise is additive when a signal is present the data have the same multivariate spherical distribution centered on the signal. We remember that this assumption generalizes to dimension larger than one the natural assumption that the noise should have a  symmetric distribution w.r. to zero i.e. negative values have the same distribution as positive ones. Exchangeability is implied by sphericity, therefore every finite set of random variables of the noise process has the same distribution of a permutation of its components. This seems a natural requirement for the noise. However sphericity implies more geometric structure. In fact e.g. the distribution  of every $n-$dimensional set of random variables of a spherical process is invariant by rotation in the $n-$dimensional Euclidean space of random variables. This too seems a natural property of the noise.  Moreover adding more noise should not modify its statistical distribution. Therefore the radial density of the multivariate spherical distribution should be an $\alpha-$stable density.

An obvious tool for solving the signal detection problem described above consists in performing a one-sample multivariate location test  $H_0$: the observed process is centered in zero, against $H_1$: the observed process is not centered in zero. It was proved in fact that it is possible to devise tests that are affine invariant and exhibit local asymptotic  optimality a la Le Cam \cite{hallin} if the joint density of the observations is elliptically symmetric and the radial density satisfies some assumptions. Unfortunately  these assumptions under the $\alpha-$stability hypothesis are not valid but in some specific cases. The idea is then to transform the original sample in order to be able to make this check.

More specifically, given an even number $n$ of complex observations, Pad\'{e} approximants of the $Z-$transform of the observed process can be computed up to order $[p-1,p]$ where $p=n/2$. Four statistics are then computed: poles, zeros, normalized residuals at the poles and normalized residuals at the zeros in the following called Pad\'{e} parameters. It turns out that all these quantities are functions of the generalized eigenvalues and eigenvectors of two pencils of random Hankel matrices. It is proved that, under $H_0$, these statistics are universal, i.e. their distribution does  not depend on the specific spherical distribution of the observations. Therefore the $\alpha-$stable radial density can be replaced by a  Gaussian one. Moreover, in the specific case of $n=2$, it is proved that the pole
 statistic satisfies the hypotheses required in  \cite{hallin}. A van der Waerden type optimal test can then be used on this parameter. A MonteCarlo experiment shows that the same test applied to the other Pad\'{e} parameters has lower power and the same is true for the same test applied to the original data and for the Hotelling test applied to the original data. Therefore it seems effortless to check the hypotheses required in  \cite{hallin} for the other Pad\'{e} parameters and the original data. We notice that the  Chernoff and Savage's result \cite[Th. 6]{hallin} comparing the van der Waerden type test and the Hotelling one on the original data does not hold in general for $\alpha-$stable data.

Finally it is proved that the poles statistic can be used to define a new optimal test a la Le Cam, whose asymptotic relative efficiency (ARE) w.r. to  the van der Waerden type test when applied to poles data is larger than one.

A MonteCarlo experiment confirms these results stressing that the advantages of the new test applied to the poles data is larger   for small values of $\alpha$ and signal-to-noise ratio (SNR).

The paper is organized as follows. In the first section the statistics are defined and their universality properties are assessed. In the second section the sphericity of the statistics is proved and the location tests are described. In the last section some simulation results are reported.

\section{Universality properties}
Let us denote random quantities by bold characters and assume that the complex-valued discrete  process
$\{\ba\}_k,\;k \in \N $, representing the signal plus white noise, is such that
 all finite sets of $\{\ba\}_k$
have an elliptically symmetric distribution. More precisely, if $\bua=[\ba_0,\dots,\ba_{n-1}]$, we assume that $\forall n = 2p,\;\;\buat=\{\Re[\bvs],\Im[\bvs]\}$ has an elliptical distribution with a density given by (see e.g. \cite{fang})
\begin{eqnarray} G(\vs;\us,\Sigma,g)=\frac{\Gamma(n)}{2 \pi^{n}}\|\vs\|_{\us,\Sigma}^{1-2 n}g(\|\vs\|_{\us,\Sigma})\label{dens}\end{eqnarray} where $g(\cdot)$ is the density of $\|\bvs\|_{\us,\Sigma}$  and $$\|\vs\|_{\us,\Sigma}=\left((\vs-\us)^T\Sigma^{-1}(\vs-\us)\right)^{1/2},\;\;\us\in\R^{2n},\;\Sigma>0\in\R^{2n\times 2n}.$$
Equivalently we can assume that
$$\buat\stackrel{d}{=}\us+\Sigma^{1/2}\bue$$
where $\us$ represents the signal and $\bue$ represents the scaled noise centered in zero with spherical distribution.

The $Z-$transform of $\{{\ba}_k\}$ is the formal random power series
$${\bf F}(z)=\sum_{k=0}^\infty \ba_k z^{-k},\;|z|>1$$
which can be extended to the unit disk by analytic continuation.
Let us denote by $[{\bf p-1,p}](z)$ the random Pad\'{e} approximant of ${\bf F}(z)$ of order $(p-1,p)$. Its poles are denoted by $\{\bxi_j\},j=1,\dots,p$ and its zeros by $\{\bzet_j\},j=1,\dots,p-1$.
The poles
can be computed  by noting that (see e.g.\cite{distrf}) they are  the generalized eigenvalues of a  pencil of square random Hankel matrices $\bU=[\bU_1,\bU_0]$  where
$$\bU_0=\left[\begin{array}{llll}
{\ba}_0 & {\ba}_{1} &\dots &{\ba}_{p-1} \\
{\ba}_{1} & {\ba}_{2} &\dots &{\ba}_{p} \\
. & . &\dots &. \\
{\ba}_{p-1} & {\ba}_{p} &\dots &{\ba}_{2p-2}
  \end{array}\right],\;
\bU_1=\left[\begin{array}{llll}
{\ba}_1 & {\ba}_{2} &\dots &{\ba}_{p} \\
{\ba}_{2} & {\ba}_{3} &\dots &{\ba}_{p+1} \\
. & . &\dots &. \\
{\ba}_{p} & {\ba}_{p+1} &\dots &{\ba}_{2p-1}\end{array}\right].$$

In \cite{zero} it was shown that  the zeros $\{\bzet_j\}$ are the poles of the random Pad\'{e} approximant of order $(p-2,p-1)$  of
$$({\bf F}(z))^{-1}=\sum_{k=0}^\infty \bbb_k z^{-k}$$
where $\{\bbb_k\}$ is defined by
$$ \bub={\bf T}^{-1}(\bua) \underline{e_1} $$
  where
$${\bf T}={\bf T}(\bua)=\left[\begin{array}{llll}
\ba_0 & 0 &\dots &0 \\
\ba_{1} & \ba_{0} &\dots &0 \\
... \\
\ba_{k-1} & \ba_{k-2} &\dots &\ba_{0}
  \end{array}\right],\;\;\;\bua=\left[\begin{array}{llll}
\ba_0  \\
\ba_{1}  \\
... \\
\ba_{k-1}
  \end{array}\right],\;\;\;\underline{e_1}=\left[\begin{array}{llll}
1  \\
0  \\
... \\
0
  \end{array}\right].$$

The zeros $\{\bzet_j\}$ are  the generalized eigenvalues of the $(p-1)\times (p-1)$ pencil ${\bf \tilde{\bU}}=[{\bf \tilde{\bU}}_1,{\bf \tilde{\bU}}_0]$ where
$${\bf \tilde{\bU}}_0=\left[\begin{array}{llll}
\bbb_2 & \bbb_{3} &\dots &\bbb_{p} \\
\bbb_{3} & \bbb_{4} &\dots &\bbb_{p+1} \\
. & . &\dots &. \\
\bbb_{p} & \bbb_{p+1} &\dots &\bbb_{2p-2}
  \end{array}\right],\;
{\bf \tilde{\bU}}_1=\left[\begin{array}{llll}
\bbb_3 & \bbb_{4} &\dots &\bbb_{p+1} \\
\bbb_{4} & \bbb_{5} &\dots &\bbb_{p+2} \\
. & . &\dots &. \\
\bbb_{p+1} & \bbb_{p+2} &\dots &\bbb_{2p-1}\end{array}\right].$$

Finally from e.g. \cite{hen2} it follows that
$${\ba}_k=\sum_{j=1}^{p}{\bf c}_j\mbox{\boldmath $\xi$}_j^k$$ therefore the residuals  ${\bf c}_j$ at the poles are given by
$$\buc=\bV(\buxi)^{-1}\bua$$ where $\bV(\buxi)$ is the random Vandermonde matrix based on $\mbox{\boldmath $\xi$}_j$. It turns out that $\bV(\buxi)^{-T}$ is the matrix of the generalized eigenvectors of  $\bU$. Analogously the residuals at the zeros are given by
$$\bud=\bV(\buzet)^{-1}\bub$$ where $\bV(\buzet)$ is the random Vandermonde matrix based on $\mbox{\boldmath $\zeta$}_j$ and $\bV(\buzet)^{-T}$ is the matrix of the generalized eigenvectors of  ${\bf \tilde{\bU}}$.

In the following we  prove that  when $\umu=\underline{0}$ and $\Sigma=I_{2n}$ (the identity matrix of order $2n$), the poles, zeros and normalized residuals $\frac{\buc}{\|\buc\|}$ an $\frac{\bud}{\|\bud\|}$ do not depend on the specific function $g(\cdot)$.
These results are  derived by the following  \cite[Theorem 2.22]{fang}

\begin{prop}
Let $\bux\in S_n^+$ where $S_n^+$ is the set of $n-$variate spherical distributions  such that $\mbox{Pr}(\bux={\bf \underline{0}})=0$. Then the distribution of a statistic $\tau(\bux)$ is invariant in $S_n^+$ provided that
$\tau(\alpha\bux)$ has the same distribution of $\tau(\bux)$ for all $\alpha>0$. In this case $\tau(\bux)$ has the same distribution of $\tau(\buy)$ where $\buy$ is an $n-$variate standard Gaussian random vector. \label{fa}
\end{prop}

\begin{theorem}
If $\umu=\underline{0}$ and $\Sigma=I_{2n}$, $\buat=\{\Re[\bvs],\Im[\bvs]\}$ is $2n$-variate spherically distributed with a density $G(\vs)=\frac{\Gamma(n)}{2 \pi^{n}}\|\vs\|^{1-2 n}g(\|\vs\|)$. Let $[{\bf p-1,p}](z)$ be the Pad\'{e} approximant of the $Z-$transform of $\{{\ba}_k\}$. Then
\begin{enumerate}
\item the marginal density of the poles and zeros of $[{\bf p-1,p}](z)$ is independent of $g(\cdot)$
and equal to the distribution obtained when $G(\vs)$ is a standard Gaussian density;
\item all statistics of normalized residuals in the poles $\frac{\buc}{\|\buc\|}$ and in the zeros $\frac{\bud}{\|\bud\|}$ are independent of $g(\cdot)$ and their distribution is equal to the distribution obtained when $G(\vs)$ is a standard Gaussian density.
\end{enumerate}
\label{th2}
\end{theorem}
\noindent \underline{proof}. Let us consider the generalized eigenvalues problem for the pencil $\bU=[\bU_1,\bU_0]$ i.e.
$$\det(\bU_1-\bxi_j\bU_0)=0.$$
By the Hankel structure of $\bU$ the solutions of this equation are invariant by multiplication of $\bua$ by a scalar $\alpha>0$. Therefore the generalized eigenvalues $\bxi_j$ of $\bU$ are statistics which satisfy the hypotheses of Proposition \ref{fa}, hence their distribution is independent of $g(\cdot)$ and equal to the distribution obtained when $G(\vs)$ is a standard Gaussian density. This concludes the proof of the first part of the first point.
\noindent Let us consider the generalized eigenvalues problem for the pencil $\tilde{\bU}=[\tilde{\bU}_1,\tilde{\bU}_0]$ i.e.
$$\det(\tilde{\bU}_1-\bzet_j\tilde{\bU}_0)=0.$$
whose solutions are invariant by multiplication of $\bub$ by a positive scalar.
But, because of the triangular Toeplitz structure of ${\bf T}$, we have
$${\bf T}^{-1}(\alpha\bua) \underline{e_1}=\alpha^{-1}{\bf T}^{-1}(\bua) \underline{e_1}=\alpha^{-1}\bub.$$
Therefore the generalized eigenvalues $\bzet_j$ of $\tilde{\bU}$ are statistics  which satisfy the hypotheses of Proposition \ref{fa}.
\noindent To prove the first part of the second point we remember that
$\buc=\bV^{-1}(\buxi)\bua$.
Let  $\tau\left(\bua\right)$ be any statistic of $\bua$ which is a function of $\frac{\buc}{\|\buc\|}$. We have
$$\tau(\bua)=\tau\left(\frac{\buc}{\|\buc\|}\right)=\tau\left(\frac{\bV^{-1}(\buxi)\bua}{\|\bV^{-1}(\buxi)\bua\|}\right)=
\tau\left(\frac{\bV^{-1}(\buxi)(\alpha\bua)}{\|\bV^{-1}(\buxi)(\alpha\bua)\|}\right).$$
But, after the first point, $\buxi$ is a function of $\bua$, invariant by multiplication for positive constants, i.e.
$\buxi(\alpha\bua)=\buxi(\bua)$. Therefore
$$\tau\left(\frac{\bV^{-1}(\buxi(\bua))(\alpha\bua)}{\|\bV^{-1}(\buxi(\bua))(\alpha\bua)\|}\right)=
\tau\left(\frac{\bV^{-1}(\buxi(\alpha\bua))(\alpha\bua)}{\|\bV^{-1}(\buxi(\alpha\bua))(\alpha\bua)\|}\right)=
\tau\left(\alpha\bua\right).$$
Therefore $\tau\left(\bua\right)$ satisfies the hypotheses of Proposition \ref{fa}.
\noindent To prove the second part of the second point we notice that if $\tau\left(\bub\right)$ is any statistic of $\bub$ which is a function of $\frac{\bud}{\|\bud\|}$ then, as before,  $\tau\left(\bub\right)=\tau\left(\alpha\bub\right)$ for all $\alpha>0$. But then if we define $\tilde{\tau}(\bua)=\bub$ we have $\tilde{\tau}(\alpha\bua)=\alpha^{-1}\bub$, and if we
  define $\hat{\tau}(\bua)=\tau(\bub)$ we get
 $$\hat{\tau}(\alpha\bua)=\tau(\tilde{\tau}(\alpha\bua)=\tau(\alpha^{-1}\bub)=\tau(\bub)=\hat{\tau}(\bua) $$
 and the thesis follows by Proposition \ref{fa}.$\;\;\;\Box$

\section{Sphericity of marginal densities of poles and residuals and location test}

The following theorems hold
\begin{theorem}
If  $\umu=\underline{0}$ and $\Sigma=I_{2n}$ and $\{\ba\}_k$ is spherical in the sense specified above, then the marginal  densities of poles and zeros of the Pad\'{e} approximant of its $Z-$transform  are spherical.
\label{th3}
\end{theorem}
\noindent \underline{proof}.
After Theorem \ref{th2} we can assume that $\{\ba\}_k$ is a complex Gaussian white noise. In \cite[Th.2]{barja} it was proved that in spherical coordinates $(\rho,\theta)$ the marginal density of a pole $\bxi=\brho e^{i\bthe}$ is a bivariate probability function $h(\rho,\theta)$ such that $h(\rho|\theta)$ does not depend on $\theta$. Therefore $\brho$ and $\bthe$ are independent and by \cite[Th.2.11]{fang} it follows that $\bxi$ has a spherical distribution.
To prove that the same property holds for the zeros we notice that the joint density of the modified process $\{\bbb\}_k$, given in \cite[Th.2]{zero}, is invariant under the transformation $$\bbb\rightarrow e^{\pm i\beta/2}\bbb .$$
Therefore the  proof of \cite[Th.2]{barja}  holds also in this case. $\;\;\Box$
\begin{theorem}
If  $\umu=\underline{0}$ and $\Sigma=I_{2n}$ and $\{\ba\}_k$ is spherical in the sense specified above, then the marginal  density of  residuals at the poles and  residuals at the zeros of the Pad\'{e} approximant of its $Z-$transform  are spherical.
\end{theorem}
\noindent \underline{proof}. We remember that, by Cramer's rule, $\bc_i=\frac{\det \bV_i}{\det \bV}$ where $\bV_i$ is the matrix obtained from $\bV$ by replacing the  $i-$th column by $\bua$. But the determinants are  measurable functions of their elements which all have spherical distribution by hypothesis and by Theorem \ref{th3}. Moreover $\det \bV$ is a.s. different from zero. Therefore $\bc_i$ is a measurable function of spherical variables. But then it is spherical by \cite[2.,pg.13]{fang}.
The same proof holds for residuals at the zeros.  $\;\;\Box$

\noindent \underline{Remark 1}. It follows by definition of sphericity that also the normalized residuals at the poles and at the zeros have a spherical distribution.

\noindent \underline{Remark 2}. As a corollary of Theorems 1 and 2 we generalize to the case of $\alpha-$stable white noise the explicit expression for the pole marginal density obtained in \cite{barja} when $n=2$ and the white noise is Gaussian.
\begin{corollary}
Let be  $n=2$ and $G(\vs;\underline{0},\sigma^2 I_4,g;\alpha)$  an $\alpha-$stable spherical density where $g(\cdot)$ is the density of $\frac{1}{\sigma}\left(\buat^T\buat\right)^{1/2}$. Then the pole  density is
$$h(z)=\frac{1}{\pi(1+|z|^2)^2},\;\;z\in\C$$
independently of $\sigma$ and $\alpha$.
\end{corollary}
\noindent \underline{proof}.
Let us assume that \cite[eq.14]{nolanm}
$$G(\vs;\underline{0},\Sigma,g;\alpha)=\frac{1}{2 \pi^{2}\sigma^4}\|\vs\|_{\underline{0},\Sigma}^{-3}g(\|\vs\|_{\underline{0},\Sigma};\alpha),\;\;r=\|\vs\|_{\underline{0},\Sigma}=
\frac{1}{\sigma}\left(\vs^T\vs\right)^{1/2},
\vs\in R^{4}$$
where by
\cite[eq.8]{nolanm} $$g(r;\alpha)=\frac{1}{2}\int_0^\infty (r t)^2 J_1(r t)e^{-(\sigma t)^\alpha} dt,$$
and $J_\nu(\cdot)$ is the Bessel function of order $\nu$.

Let us consider   the change of variables $(a_0,a_1)\rightarrow(\gamma,\lambda)$ given by
$a_0=\gamma,\;a_1=\gamma\lambda$ with real Jacobian $|\gamma|^2$ and $r=\frac{1}{\sigma}\left(\vs^T\vs\right)^{1/2}=\frac{|\gamma|}{\sigma}(1+|\lambda|^2)^{1/2}$.
But then
the pole marginal density is given by
$$E\left[\delta\left(z-\frac{a_1}{a_0}\right)\right]=\int_{\R^4}|\gamma|^2\delta(z-\lambda)
G([\Re{\gamma},\Im{\gamma},\Re{\lambda},\Im{\lambda}];\alpha)d\Re{\gamma}d\Im{\gamma}d\Re{\lambda}d\Im{\lambda}=$$
$$\int_{\R^2}|\gamma|^2
G([\Re{\gamma},\Im{\gamma},\Re{z},\Im{z}];\alpha)d\Re{\gamma}d\Im{\gamma}=$$ $$\frac{1}{2 \pi^{2}\sigma}\int_{\R^2}\frac{1}{|\gamma|(1+|z|^2)^{3/2}}
g\left(\frac{1}{\sigma}|\gamma|(1+|z|^2)^{1/2};\alpha\right)d\Re{\gamma}d\Im{\gamma}.$$
Let us consider the change of variables $(\Re{\gamma},\Im{\gamma})\rightarrow(\rho,\theta)$ given by
$$\Re{\gamma}= \frac{\sigma \rho}{\sqrt{1+|z|^2}}\cos{\theta},\;\;\Im{\gamma}= \frac{\sigma \rho}{\sqrt{(1+|z|^2}}\sin{\theta}  $$ with Jacobian $\frac{\sigma^2\rho}{1+|z|^2}$. Let be  $K_1=(1+|z|^2)^{3/2},\;\;K_2= \frac{1}{\sigma}(1+|z|^2)^{1/2}$, then we get
$$\frac{1}{2 \pi^{2}\sigma}\int_{\R^2}\frac{1}{|\gamma|(1+|z|^2)^{3/2}}
g\left(\frac{1}{\sigma}|\gamma|(1+|z|^2)^{1/2};\alpha\right)d\Re{\gamma}d\Im{\gamma}=$$
$$\frac{1}{2 \pi^{2}\sigma} \int_{\R^2}\frac{\rho}{K_2^2}\frac{K_2}{\rho K_1 }
g\left(\rho;\alpha\right)d\rho d\theta=\frac{1}{2 \pi^{2}\sigma} \frac{1}{K_2}\frac{1}{ K_1 }\int_{\R^2}
g\left(\rho;\alpha\right)d\rho d\theta=$$ $$\frac{1}{ \pi\sigma} \frac{1}{K_2}\frac{1}{ K_1 }\int_{\R^2}
g\left(\rho;\alpha\right)d\rho =\frac{1}{ \pi\sigma} \frac{1}{K_1K_2}=\frac{1}{\pi(1+|z|^2)^2}.\;\;\;\Box$$

We now consider the case $n=2$ and   the pole statistic because it is the most promising one when $\alpha$ and the SNR are small, as it will be shown in the following. We want to show that for this statistic the conditions of applicability of the van der Waerden type optimal test, proposed in \cite{hallin}, hold.

Let us define the radial function $f(r)$ of an elliptic density $G(\vs;\umu,\Sigma,g),\;\vs\in\R^n$ as the function which satisfies the equation
$$g(r)=\frac{r^{n-1}}{\nu_{n-1}}f(r)$$
where $g(r)$ is the density of $\|\bvs\|_{\umu,\Sigma}$ and
$$\nu_k=\int_0^\infty r^k f(r) dr.$$
If $n=2$ and the complex data are $\ba_k=c\xi^{k-1}+\sigma\bue_k,\;k=1,2;$ where  $\xi=(\xi_R,\xi_I)=\frac{s_1}{s_0}$ denotes the true pole,  and $\rho=\frac{|c|^2}{\sigma^2}$ denotes the $SNR$ , then in the limits $\rho\rightarrow 0,\;\;|\xi|\rightarrow 1$ the pole statistic satisfies the hypotheses (A1) and (A2') in  \cite{hallin}. In fact
\begin{theorem}
If $n=2$ and $G(\vs;\underline{s},\frac{\sigma^2}{2} I_4,g;\alpha)$  an $\alpha-$stable spherical density where $g(\cdot)$ is the density of $\frac{\sqrt{2}}{\sigma}\left(\buat^T\buat\right)^{1/2}$  then the pole density can be approximated, for $\rho\rightarrow 0,\;\;|\xi|\rightarrow 1$, by the spherical density
$$\tilde{h}(z,\rho)= \frac{1}{\pi  \left(1+|z- \rho\Gamma(1+\frac{2}{\alpha})\xi |^2\right)^2},  $$
the pole radial function is $f(r)=\frac{1}{(1+r^2)^2},$ its first moment $\nu_1$ is finite, and
$(f^{1/2})'\in L^2(\R^+_0,\nu_1)$ where $L^2(\R^+_0,\nu_1)$ is the space of square-integrable function w.r. to the Lebesgue measure with weight $r$.
\end{theorem}
\noindent \underline{proof}.
Let us assume that \cite[eq.14]{nolanm}
$$G(\vs;\us,\Sigma,g;\alpha)=\frac{2}{ \pi^{2}\sigma^4}\|\vs\|_{\us,\Sigma}^{-3}g(\|\vs\|_{\us,\Sigma};\alpha),\;\;r=\|\vs\|_{\us,\Sigma}=
\frac{\sqrt{2}}{\sigma}\left((\vs-\us)^T(\vs-\us)\right)^{1/2},
\vs\in R^{4}$$
where by
\cite[eq.8]{nolanm} $$g(r;\alpha)=\frac{1}{2}\int_0^\infty (r t)^2 J_1(r t)e^{-\left(\frac{\sigma}{\sqrt{2}} t\right)^\alpha} dt,$$
and $J_m(\cdot)$ is the Bessel function of order $m$.

Let us consider   the change of variables $(a_0,a_1)\rightarrow(\gamma,\lambda)$ given by
$a_0=\gamma,\;a_1=\gamma\lambda$ with real Jacobian $|\gamma|^2$ and $$r=\|\vs\|_{\us,\Sigma}=\left((\gamma-\mu)^H\frac{2(1+|\lambda|^2)}{\sigma^2}(\gamma-\mu)+q\right)^{1/2}$$
where (\cite[Lemma 1.1]{qam})
$$\mu=\frac{s_0+\overline{\lambda}s_1}{1+|\lambda|^2},\;\;q=|s_0|^2+|s_1|^2-\frac{|s_0+\overline{\lambda}s_1|^2}{1+|\lambda|^2}=
|s_0|^2\frac{|\lambda-\xi|^2}{1+|\lambda|^2}$$
But then
the pole marginal density is given by
$$h_2(z;\alpha)=E\left[\delta\left(z-\frac{a_1}{a_0}\right)\right]=\int_{\R^4}|\gamma|^2\delta(z-\lambda)
G([\Re{\gamma},\Im{\gamma},\Re{\lambda},\Im{\lambda}];\alpha)d\Re{\gamma}d\Im{\gamma}d\Re{\lambda}d\Im{\lambda}=$$
$$=\int_{\R^2}|\gamma|^2
G([\Re{\gamma},\Im{\gamma},\Re{z},\Im{z}];\alpha)d\Re{\gamma}d\Im{\gamma}$$
$$=\frac{2}{ \pi^{2}\sigma^4 }\int_{\R^2}
\frac{g\left(r;\alpha\right)}{r^3}d\Re{\gamma}d\Im{\gamma}.$$
Let us substitute $\lambda$ with $z$ in $\mu$ and $q$ and let be $a=\sqrt{\frac{2q}{\sigma^2}}$. Let us consider the change of variables $(\Re{\gamma},\Im{\gamma})\rightarrow(r,\theta)$ given by
$$\Re{\gamma}= \mu_R + \frac{\sigma}{\sqrt{2(1 +|z|^2)}} \sqrt{r^2 - a^2} \cos{\theta} $$
$$\Im{\gamma}= \mu_I + \frac{\sigma}{\sqrt{2(1 +|z|^2)}} \sqrt{r^2 - a^2} \sin{\theta} $$
 with Jacobian $\frac{\sigma^2 r}{2(1+|z|^2)}$. Then we get
$$h_2(z;\alpha)=\frac{1}{\pi^2\sigma^2}\frac{1}{(1 + |z|^2)} \int_{a}^{\infty}\int_{-\pi}^{\pi}\frac{g(r;\alpha)}{r^2}  (\Re{\gamma}^2 +\Im{\gamma}^2) d\theta dr$$
$$
=\frac{2}{\pi \sigma^2} \frac{ \left( \left(1+|z|^2\right) |\mu|^2- q\right)}{\left(1+|z|^2\right)^2 }\int_a^{\infty }\frac{ g(r;\alpha )}{r^2} dr +
\frac{1}{\pi \left(1+|z|^2\right)^2} \int_a^{\infty } g(r;\alpha) dr$$
$$=\frac{2}{\pi}\rho\frac{|1+\overline{z}\xi|^2-|z-\xi|^2}{(1 + |z|^2)^3}\int_a^{\infty }\frac{ g(r;\alpha )}{r^2} dr+
\frac{1}{\pi \left(1+|z|^2\right)^2} \int_a^{\infty } g(r;\alpha) dr.$$
Let us define
$$K_1=\frac{|1+\overline{z}\xi|^2-|z-\xi|^2}{(1 + |z|^2)^3},\;\;K_2=\frac{1}{\pi \left(1+|z|^2\right)^2},\;\;K_3=
\frac{2|z-\xi|^2}{1 + |z|^2}.$$
We get
$$h_2(z,\rho;\alpha)=\frac{2}{\pi}\rho K_1\int_{\sqrt{K_3\rho}}^{\infty }\frac{ g(r;\alpha )}{r^2} dr+K_2\int_{\sqrt{K_3\rho}}^{\infty } g(r;\alpha) dr.$$
In order to compute the first order Taylor series of $h_2(z,\rho;\alpha)$ around $\rho=0$ we compute its derivative
$$\frac{\partial h_2(z,\rho;\alpha)}{\partial \rho}=\frac{-(2K_1+K_2K_3\pi)g(\sqrt{K_3\rho};\alpha )+4K_1\sqrt{K_3\rho}
\int_{\sqrt{K_3\rho}}^{\infty }\frac{ g(y;\alpha )}{y^2}dy}{2 \pi\sqrt{K_3\rho}}.$$
But then
$$\lim_{\rho\rightarrow 0}\frac{\partial h_2(z,\rho;\alpha)}{\partial \rho}=-\frac{(2K_1+K_2K_3\pi)}{2\pi}\lim_{\rho\rightarrow 0}
\frac{g(\sqrt{K_3\rho};\alpha )}{\sqrt{K_3\rho}}+\frac{4K_1}{2\pi}\int_0^\infty\frac{ g(y;\alpha )}{y^2}dy.$$
From \cite[eq.12]{nolanm}
$$\lim_{\rho\rightarrow 0}
\frac{g(\sqrt{K_3\rho};\alpha )}{\sqrt{K_3\rho}}=\lim_{\rho\rightarrow 0}K_3\rho\frac{g(\sqrt{K_3\rho};\alpha )}{(\sqrt{K_3\rho})^3}=
\lim_{\rho\rightarrow 0}K_3\rho \frac{4\Gamma(4/\alpha)}{\alpha}=0$$
and, from \cite[eq.24]{nolanm},
$$\int_0^\infty\frac{ g(y;\alpha )}{y^2}dy=\frac{1}{2}\Gamma(1+2/\alpha)$$
therefore
$$\lim_{\rho\rightarrow 0}\frac{\partial h_2(z,\rho;\alpha)}{\partial \rho}=K_1\frac{\Gamma(1+2/\alpha)}{\pi}.$$
But then, remembering that  $g(\cdot;\alpha )$ is a density,
$$h_2(z,\rho;\alpha)=K_2+\rho \frac{\Gamma(1+2/\alpha)}{\pi}K_1+O\left(\rho ^2\right)$$ $$=
\frac{1}{\pi \left(1+|z|^2\right)^2}+\rho \Gamma(1+2/\alpha)\frac{|1+\overline{z}\xi|^2-|z-\xi|^2}{\pi(1 + |z|^2)^3}+O\left(\rho ^2\right)$$
$$=\frac{1}{\pi \left(1+|z|^2\right)^2}+\rho \Gamma(1+2/\alpha)\left(\frac{ \left(|\xi|^2-1\right) \left(|z|^2-1\right)}{\pi
   \left(1+|z|^2\right)^3}+\frac{4   \Re[\overline{z}\xi]}{\pi\left(1+|z|^2\right)^3}\right)+O\left(\rho ^2\right)$$
Let us consider the elliptic density
$$\tilde{h}(z,\rho)= \frac{1}{\pi  \left(1+|z- \rho\Gamma(1+2/\alpha)\xi |^2\right)^2}$$
and    its first order Taylor series  around $\rho=0$
$$\tilde{h}(z,\rho)=\frac{1}{\pi  \left(1+|z|^2\right)^2}+\rho\Gamma(1+2/\alpha)\frac{4   \Re[\overline{z}\xi]}{\pi\left(1+|z|^2\right)^3}+O\left(\rho ^2\right) . $$
 The pole density is then well approximated by a spherical density centered in $\rho\Gamma(1+2/\alpha)\xi$  when $\rho\rightarrow 0,\;\;|\xi|\rightarrow 1$.
  But then if $f(r)=\frac{1}{(1+r^2)^2}$ we have $$\nu_1=\int_0^\infty r f(r) dr = \frac{1}{2}.$$ Hence $f(r)$ is the pole radial function. Moreover
$$(f^{1/2})'=-g(r)$$ and
$$\int_0^\infty[(f^{1/2})']^2 r dr= \int_0^\infty \frac{4 r^3}{\left(r^2+1\right)^4}  dr=\frac{1}{3}.\;\;\;\Box$$

As a consequence of this theorem  the hypotheses (A1) and (A2')  in  \cite{hallin} are satisfied and the LAN property required in \cite[Prop.2]{hallin} holds and this is enough to apply the theory developed there to the pole statistic.

Let us denote by ${\mathcal H^{(m)}}(\umu,\sigma^2I_2,f)$ the hypothesis under which the observations have joint density $\prod_{i=1}^m G(\vs_i;\umu,\sigma^2I_2,f)$ and by
$${\mathcal H_0^{(m)}}(\umu,\Sigma,f)= \bigcup_\Sigma\bigcup_f{\mathcal H^{(m)}}(\umu,\Sigma,f)$$
$${\mathcal H_1^{(m)}}(\umu,\Sigma,f)= \bigcup_{\umu\ne \underline{0}}\bigcup_\Sigma\bigcup_f{\mathcal H^{(m)}}(\umu,\Sigma,f)$$
the testing problem we are interested in.
We want to show now that the pole  radial function can be used as a score function giving rise to a new optimal test in the class proposed in \cite[eq. (5)]{hallin}. More specifically
let us define a bivariate sample of dimension $m$  from an elliptical density with radial function $f(r)=\frac{1}{(1+r^2)^2}$ by $\{\bux_1,\dots,\bux_m\}$ and let us consider the test statistic  defined by
$$Q^{(m)}_{f_*}(\underline{0})=\frac{2}{m E[J^2_{2,f_*}(\bu)]}\sum_{i,j=1}^m J_{2,f_*}\left(\frac{\hat{R}_i(\underline{0})}{m+1}\right)
J_{2,f_*}\left(\frac{\hat{R}_j(\underline{0})}{m+1}\right)\cos(\pi p_{ij}^{(m)}(\underline{0}))$$
where $\bu$ is uniform in $(0,1)$; $J_{2,f_*}=-2\frac{(f_*^{1/2})'}{(f_*^{1/2})'} \circ \tilde{g}^{-1}$, $\circ$ denotes functions composition and $\tilde{g}(r)$ is the distribution function associated to $g(r)$; $p_{ij}^{(m)}(\underline{0})$ are the normalized interdirections \cite{rand}; $\hat{R}_j(\underline{0})$ are the pseudo-Mahalanobis ranks of the sample $\{\bux_1,\dots,\bux_m\}$.
We have
\begin{lemma}
If $f_*(r)=\frac{1}{(1+r^2)^2}$ then $J_{2,f_*} (u)=4\sqrt{u(1-u)}$ and
$$\int_0^1{|J_{2,f_*} (u)|^{2+\delta}}du=\frac{\sqrt{\pi } 2^{\delta +1} \Gamma \left(\frac{\delta }{2}+2\right)}{\Gamma
   \left(\frac{\delta +5}{2}\right)}
<\infty,\;\forall \delta>0$$
\end{lemma}
therefore assumption (A3) holds true and, by using an affine-invariant scatter estimator as defined e.g. in \cite{tyler}, also assumption (A4) is satisfied. Therefore Prop. 3 and Prop. 4 in \cite{hallin} are true i.e.

\begin{theorem}(Hallin, Paindaveine \cite[Prop. 4]{hallin})
  The sequence of tests ${\mathcal P}^{(m)}$ rejecting the null hypotesis $H_0$ whenever $Q^{(m)}_{f_*}(\underline{0})$ exceeds the $(1-\beta)$ quantile $\chi^2_{2,(1-\beta)}$ of a chi-square distribution with $2$ degrees of freedom and $f_*(r)=\frac{1}{(1+r^2)^2}$

  (i) has asymptotic level $\beta$

  (ii) is locally asymptotically maxmin, at asymptotic level $\beta$, for
  ${\mathcal H_0^{(m)}}(\umu,\Sigma,f)$ against alternatives of the form
$$\bigcup_{\umu\ne \underline{0}}\bigcup_\Sigma{\mathcal H^{(m)}}(\umu,\Sigma,f_*).$$

 \end{theorem}

Let us denote by ${\mathcal W}^{(m)}$  the sequence of tests of the van der Waerden type defined by $Q^{(m)}_{f_*}(\underline{0})$ when $f_*(r)=e^{-\frac{r^2}{2}}$.
If $ARE ({\mathcal P}^{(m)},{\mathcal W}^{(m)})$ denotes the asymptotic Pitman relative efficiency (see e.g. \cite[sec. 14.4]{vaart} of ${\mathcal P}^{(m)}$ w.r. to ${\mathcal W}^{(m)}$, we have

\begin{theorem}
 If $\{\bux_1,\dots,\bux_m\}$ is a bivariate sample from the pole distribution when $n=2$ then $$ARE ({\mathcal P}^{(m)},{\mathcal W}^{(m)})\approx 1.44$$
\end{theorem}

\noindent \underline{proof}.
We have
$$ARE ({\mathcal P}^{(m)},{\mathcal W}^{(m)})=\frac{C^2_{{\mathcal P}^{(m)}}(f_*,f)}{C_{{\mathcal P}^{(m)}(f_*,f_*)}}\cdot \frac{C_{{\mathcal W}^{(m)}(f_*,f_*)}}{C^2_{{\mathcal W}^{(m)}(f_*,f)}}$$
where
$$C(f_1,f_2)=\int_0^1 J_{2,f_1} (u)J_{2,f_2} (u) du.$$
But if $f_*(u)=\frac{1}{(1+r^2)^2}=f(u)$ then
$$C_{{\mathcal P}^{(m)}(f_*,f_*)}=C_{{\mathcal P}^{(m)}(f_*,f)} =16\int_0^1 {u(1-u)} du= \frac{8}{3}$$
and if $f_*(u)=e^{-\frac{r^2}{2}}$ and $f(u)=\frac{1}{(1+r^2)^2}$ then
$$C_{{\mathcal W}^{(m)}(f_*,f_*)}= 2\int_0^1 \log\left(\frac{1}{1-u}\right) du=2  $$
and
$$C_{{\mathcal W}^{(m)}(f_*,f)}=4\int_0^1 \sqrt{2\log\left(\frac{1}{1-u}\right)} \sqrt{u(1-u)}du \approx 1.92.\;\;\;\Box$$

 As a final remark we notice  that the  Chernoff and Savage's result \cite[Th. 6]{hallin} comparing the van der Waerden type test and the Hotelling one on the original data does not hold in general for $\alpha-$stable data. In fact assumption (A1') of \cite{hallin} does not hold because (see e.g. \cite[ch.VI]{feller2})
\begin{lemma}
If $f(r)$ is an $\alpha-$stable radial function with $\alpha<2$ and $k>0$ then $\int_0^\infty r^{k+1} f(r) dr =\infty.$
\end{lemma}

\section{Simulations}

In this section some of the claims done in the previous sections are checked by MonteCarlo simulations. In fig.1, ${\mathcal W}^{(m)}$, applied to pole data, has larger power than  the same test applied to the other Pad\'{e} parameters. This justify the use of the pole statistic.  Moreover ${\mathcal W}^{(m)}$, applied to pole data, performs better  than applied to the original data at least when $\alpha<0.5$ and the smaller $\alpha$ the larger the advantage to use the van der Waerden type test on pole data (figs.2-4).
Finally  in figs.2-4 is shown that the power of the test ${\mathcal P}^{(m)}$ applied to pole data is larger than those of the other tests considered above when $\alpha$ and the SNR are small.
The theoretical power of the Hotelling $T^2$ test applied to Gaussian data is also reported for comparison in figs.2-4.

\newpage

\begin{figure}
\begin{center}
\hspace{1.7cm}{\fbox{\psfig{file=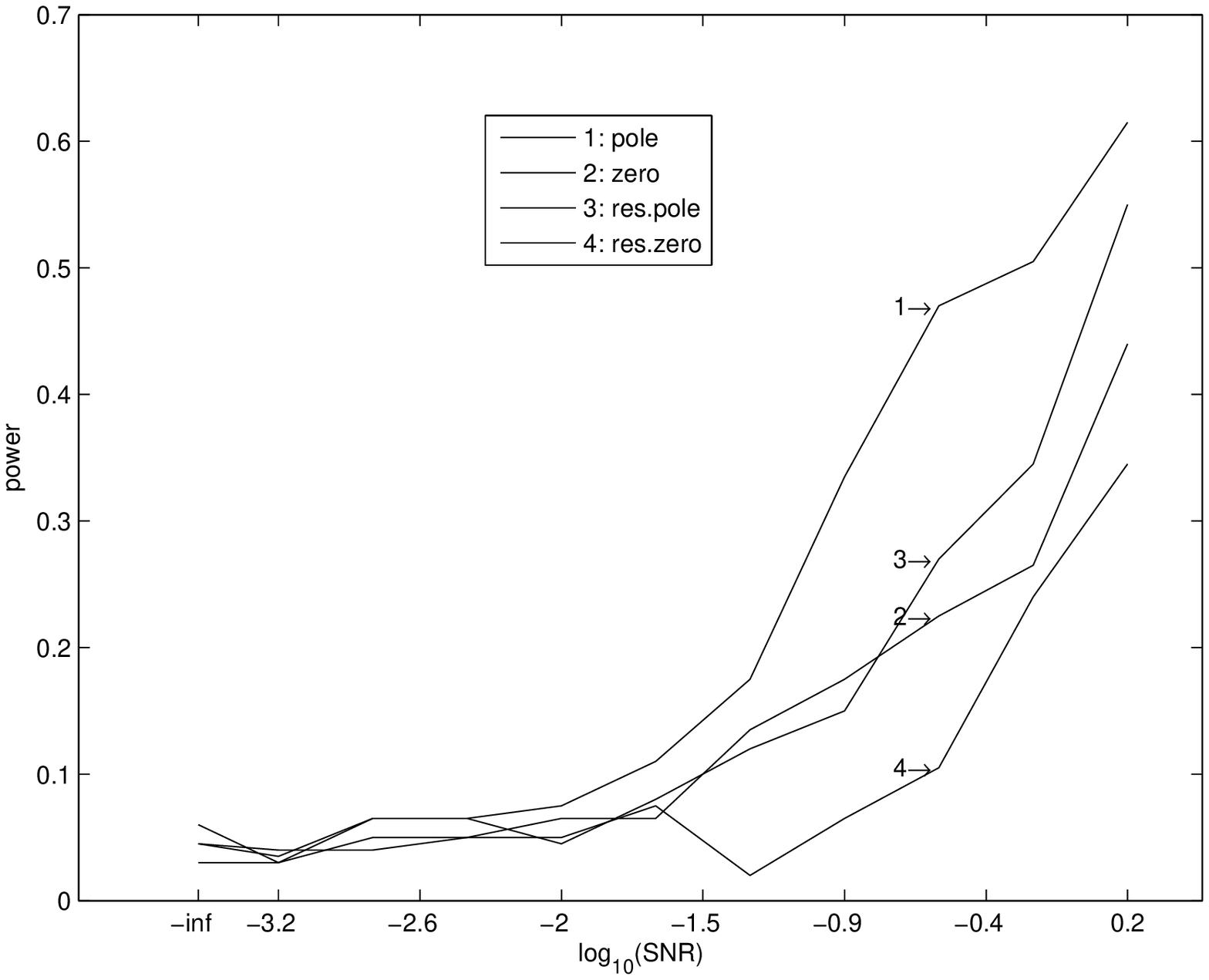,height=16cm,width=16cm}}}
\end{center}
\caption{Power of ${\mathcal W}^{(m)}$ applied to Pad\'{e} parameters as a function of  the $SNR$ when $\alpha = 0.1$, the   significance level is $\beta=0.05$, $m=100$ and the number of Montecarlo samples is $200$. } \label{fig1}
\end{figure}

\begin{figure}
\begin{center}
\hspace{1.7cm}{\fbox{\psfig{file=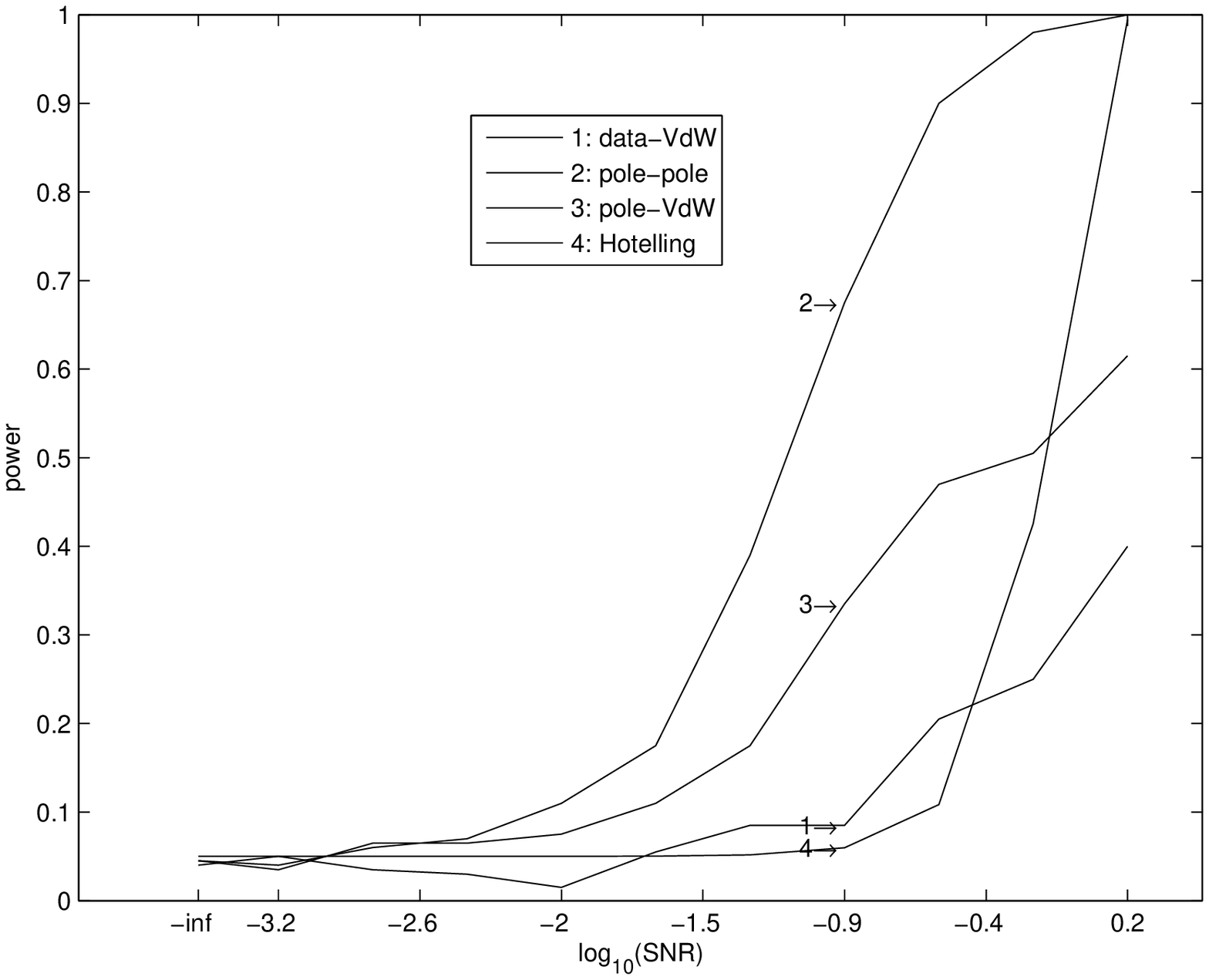,height=16cm,width=16cm}}}
\end{center}
\caption{Power of ${\mathcal W}^{(m)}$ and ${\mathcal P}^{(m)}$ tests applied to the original and pole data. The theoretical power of the Hotelling test applied to Gaussian data is also shown. The powers are reported  as  functions of the  $SNR$ when $\alpha = 0.1$, the  significance level is $\beta=0.05$, $m=100$ and the number of Montecarlo samples is $200$.} \label{fig2}
\end{figure}

\begin{figure}
\begin{center}
\hspace{1.7cm}{\fbox{\psfig{file=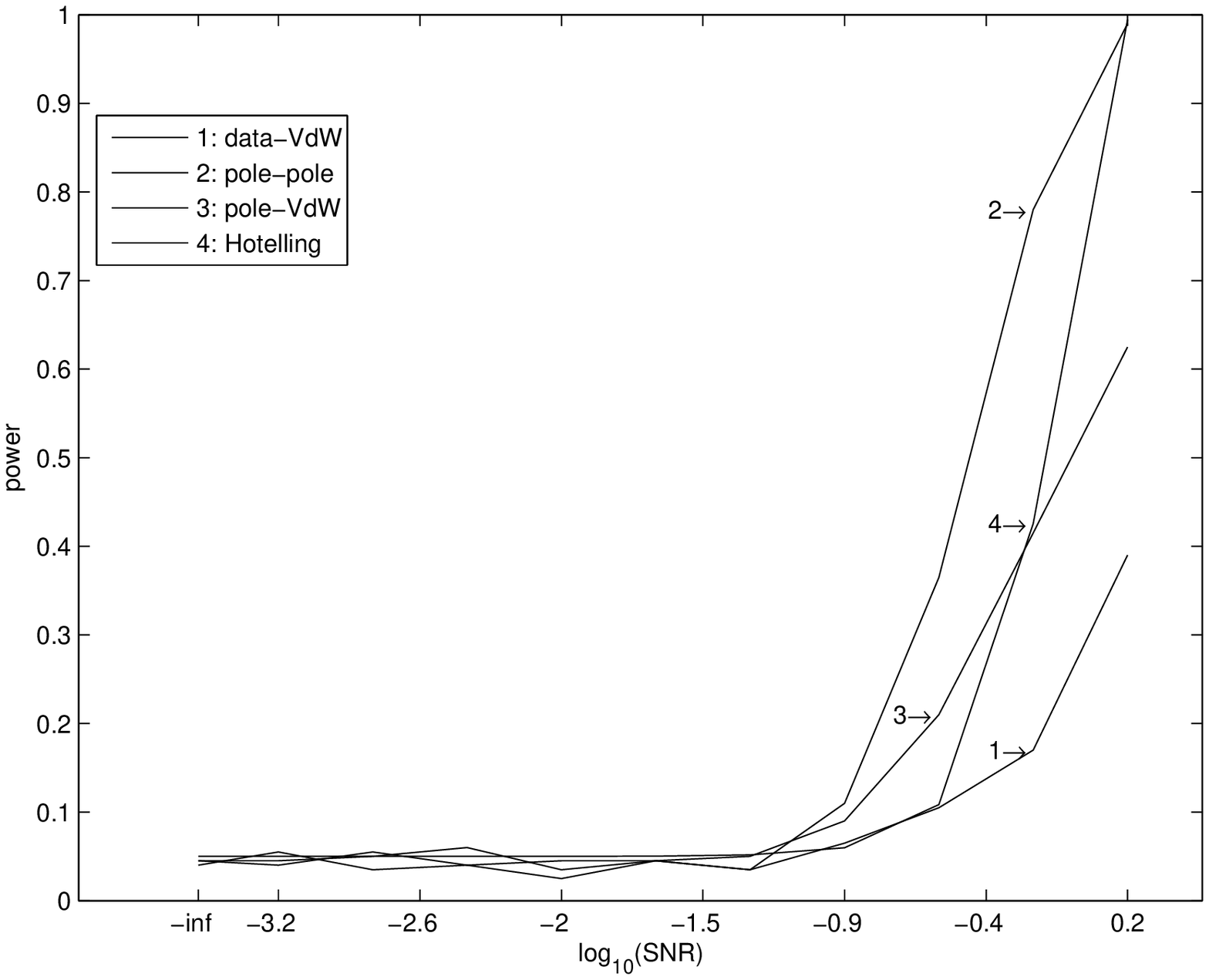,height=16cm,width=16cm}}}
\end{center}
\caption{Power of ${\mathcal W}^{(m)}$ and ${\mathcal P}^{(m)}$ tests applied to the original and pole data. The theoretical power of the Hotelling test applied to Gaussian data is also shown. The powers are reported  as  functions of  the $SNR$ when $\alpha = 0.2$, the  significance level is $\beta=0.05$, $m=100$ and the number of Montecarlo samples is $200$.} \label{fig3}
\end{figure}

\begin{figure}
\begin{center}
\hspace{1.7cm}{\fbox{\psfig{file=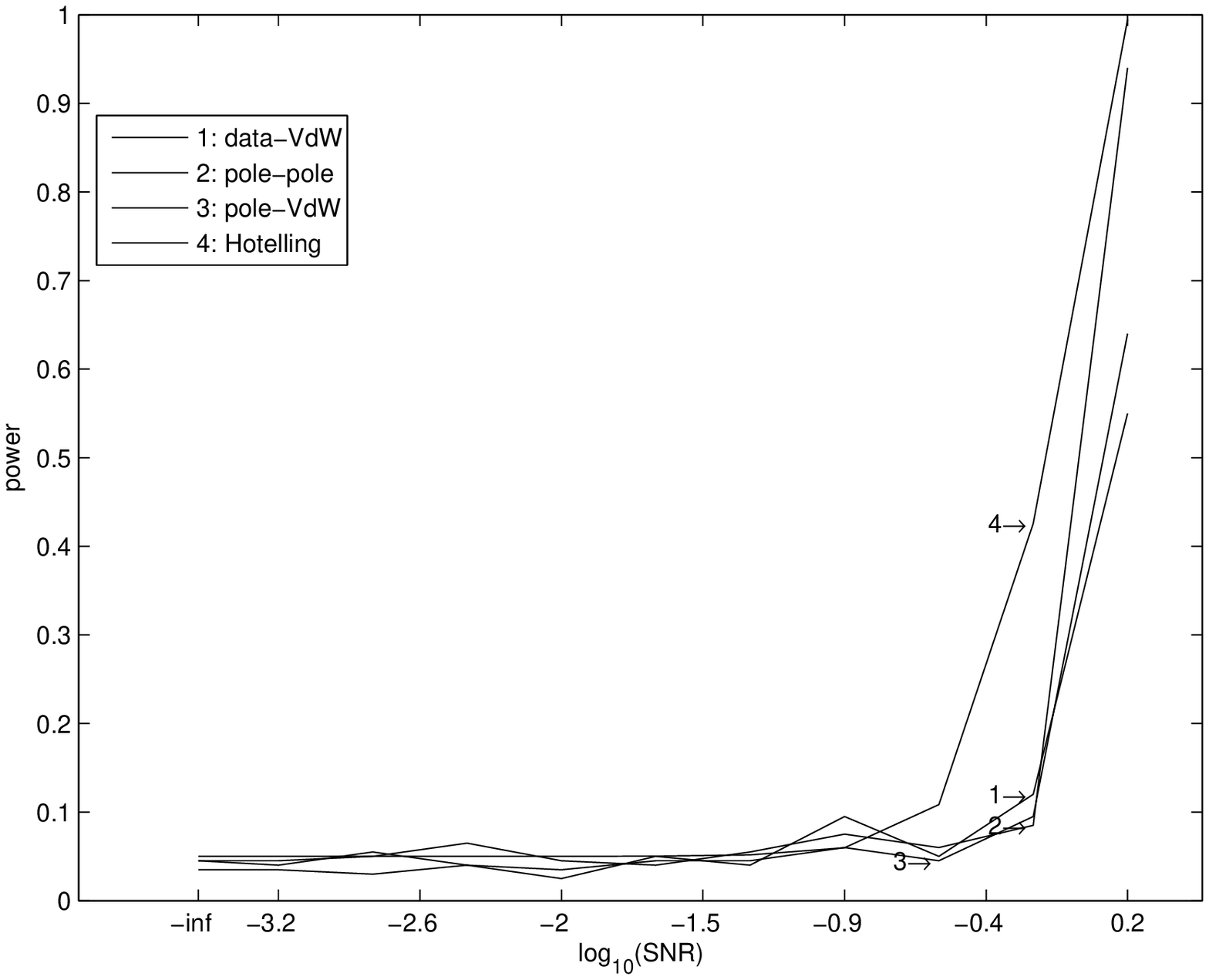,height=16cm,width=16cm}}}
\end{center}
\caption{Power of ${\mathcal W}^{(m)}$ and ${\mathcal P}^{(m)}$ tests applied to the original and pole data. The theoretical power of the Hotelling test applied to Gaussian data is also shown. The powers are reported  as  functions of the  $SNR$ when $\alpha = 0.5$, the  significance level is $\beta=0.05$, $m=100$ and the number of Montecarlo samples is $200$.} \label{fig4}
\end{figure}

\end{document}